\documentstyle[12pt]{article}
\begin{document}

\title{Approximate Rational Arithmetics and Arbitrary Precision
Computations\footnote{The work was supported by the RFBR grant
99--01--01198.\newline
Russian version: {\it Computational Technologies\/} {\bf 5}:5 (2000).}}
\author{G.~L.~Litvinov, A.~Ya.~Rodionov, A.~V.~Chourkin}
\date{}

\maketitle

\centerline{}

\centerline{\bf Abstract}

\medskip
{\it  We describe an approximate rational arithmetic with round-off errors
(both absolute and relative) controlled by the user. The rounding procedure
is based on the continued fraction expansion of real numbers. Results
of computer experiments are given in order to compare efficiency and
accuracy of different types of approximate arithmetics and rounding
procedures.}

\medskip


\medskip

{\bf 1. Introduction.} Problems of validity and reliability of
calculations (including the analysis of round-off errors) are becoming
more and more important recently, partly due to the steady growth of
computer power. Roughly speaking, the main disadvantage of the standard
floating point arithmetic is, that relative round-off error only can be
controlled during calculations. In some cases (for example, in summation
of series and subtraction of nearly equal numbers) this disadvantage can
lead to a loss of accuracy and even to absolutely incorrect results. So,
if the result of calculations depends on the errors in input data and
round-off errors critically (for example, in the case of solving ill-posed
equations, study of stability of solutions etc.), then it is reasonable
to use calculations with multiple and even arbitrary precision.

An appealing way to improve the accuracy of calculations is to use
different versions of rational arithmetics, which work with rational
numbers of the form $\frac{p}{q}$, where $p$ and $q$ are integer numbers
($q>0$). It is possible to use the exact rational arithmetic (see, e.g.
\cite{1}), but, as a rule, it leads to an explosive growth of both
calculation time and storage space since magnitudes (and lengths) of
numerators and denominators of computed numbers grow very fast.
The approximate rational arithmetics with fixed slash (maximum of
lengths of the numerator and the denominator is fixed) or floating
slash (sum of these lengths is fixed) were investigated in detail
earlier (see \cite{2, 3}). These rounding procedures use the representation
of continued fractions. There also exist rougher rounding procedures,
which use only a fixed number of top digits (the other digits are replaced
by zeros) but it sometimes leads to relatively large rounding errors.

Here we suggest a new modification of approximate rational arithmetic
with a more natural and accurate rounding procedure. The user defines
values $\Delta$ and $\delta$ of absolute and relative error such that
$0\leq\Delta\leq\infty$, $0\leq\delta\leq\infty$.  In particular, in
the case $\Delta=\delta=0$ we obtain the exact rational arithmetic.
If $\delta=\infty$, then only absolute rounding error $\Delta$ is fixed.
If only relative error  $\delta$ is fixed, then we obtain approximately
the same picture as for the floating-point arithmetic. This rounding procedure
is applied to a fraction if lengths of its numerator and denominator
exceed a number $M$ specified by the user. In this case the initial
rational number is replaced by its best approximation in the form of a
convergent of a continued fraction within given errors $\Delta$ and $\delta$.
The result of rounding is always an uncancellable fraction, and sometimes
it can coincide with the initial number.

This type of arithmetic was originally implemented by means of the REDUCE
computer algebra system  and  was used for constructing arbitrary rational
approximations to functions of one variable \cite{4}. In this paper the
analysis of accuracy and efficiency of different modifications of
approximate rational arithmetics is based on computer experiments, which
are implemented by means of the $C^{++}$ language (in an object oriented
form within the framework of the project outlined in \cite{5}). Note that
Yu.~V.~Matijasevich suggested to apply an approximate rational
arithmetic of arbitrary precision for his {\it a posteriori interval
analysis } \cite{6,7}.

Below we describe an algorithm of rounding and construction of best
approximations (using the convergents of continued fractions). This
method is compared with other methods by means of a computer experiment.
We give also an estimation of the number of components of continued
fractions in dependence on the accuracy of rounding. In particular, we
show that  increase in accuracy of calculations does not lead to an
explosive growth of calculation time. The main result of computer
experiments is that this algorithm of rounding provides significantly
higher accuracy of calculations in comparison with other modifications
of rational arithmetics  that require  comparable time.

{\bf 2. Continued fractions and the approximate rational arithmetic.}
Recall some basic notions of the theory of continued fractions. Denote
by $[ a_0;a_1,a_2,...] $ a continued fraction of the form
$$
a_0+\frac{1}{a_1+\frac{1}{a_2+...}}.\eqno(1)
$$
Any non-negative rational number $\frac{p}{q}$ has a unique canonical
representation in the form of a finite continued fraction
$$
\frac{p}{q}=[ a_0;a_1,a_2,...,a_n],\eqno(2)
$$
where all $a_i$ ($i=0,...,n$) are non-negative integers,
$a_i\ge1$ if  $i=1$, \dots , $n-1$ and $a_n\ge2$ if  $n\ge1.$ Irrational
numbers can also be represented in the form (1) as infinite continued
fractions. {\it A convergent} $\frac{p_k}{q_k}$ {\it of the order k of a
continued fraction}  is defined for the decomposition (2) by the equality
$$
\frac{p_k}{q_k}=[ a_0;a_1,a_2,...,a_k],\eqno(3)
$$
where $k\le n$. It is clear from (3) that the convergent of a continued
fraction $\frac{p_n}{q_n}$ coincides with $\frac{p}{q}$. From the theory
of continued fractions \cite{8} it is well known that the convergent of
a continued fraction (3) is a best approximant to the number (2) in the
following sense: for any fraction $\frac{r}{s}$ such that $0<s\le q_k$
and $\frac{r}{s}\ne\frac{p_k}{q_k}$, the following inequality holds:
$$
\Big|\frac{r}{s}-\frac{p}{q} \Big| > \Big|\frac{p_k}{q_k}-\frac{p}{q}\Big|.
$$
The only (trivial) counterexample is $\frac{p}{q}=a_0+\frac{1}{2}$; in
this case $a_0$ and $a_0+1$ approximate the number $p/q$ equally well.
For any convergent $\frac{p_k}{q_k}$ in the case of $k\ne0$ and $k<n$,
the following inequality holds:
$$
\frac{1}{q_k(q_k+q_{k+1})}<\Big|\frac{p}{q}-\frac{p_k}{q_k}\Big| \le
\frac{1}{q_kq_{k+1}} \eqno(4)
$$
A similar property is correct for infinite continued fractions of the
form (1) corresponding to irrational numbers. Using the inequality (4)
we obtain an efficient algorithm for estimation of the rounding error.
These properties of convergents make continued fractions an ideal tool
for construction of approximate rational arithmetic.

Recall the following algorithm for constructing of a convergents of a
continued fraction \cite{8}. Let the initial fraction be $\frac{p}{q}$,
where $p, q$ are integer numbers, $p>0$, $q>1$ (if $p/q<0$, then we work
with the fraction $\vert p/q\vert$ and then multiply result by $-1$).
\begin{itemize}
\item{} Let the initial condition be defined by $b_{-2}=p$, $ b_{-1}=q$,
$ p_{-2}=0$, $ p_{-1}=1$, $ q_{-2}=1$, $q_{-1}=0.$
\item{}  For $i=1,2,...$, the values of $a_i$ and $b_i$ are consecutively
computed as the quotient and the remainder obtained when $b_{i-2}$ is
divided by $b_{i-1}$ respectively:
$$ b_{i-2}=a_ib_{i-1}+b_i .$$
 \item{} The numerator and the denominator of the convergent of order $i$
of the continued fraction are given in the recurrent form:
$$
p_i=a_ip_{i-1}+p_{i-2},
$$
$$
q_i=a_iq_{i-1}+q_{i-2}.
$$
\item{} If $b_i=0$, then the convergent of the continued fraction coincides
with the initial fraction $\frac{p}{q}$ and the procedure terminates.
\item{} At each step (at each {\it i=0,1,..}), a criterion of accuracy (see
below) is checked, and if the result satisfies the criterion of accuracy
then the procedure terminates, otherwise we  perform the next step with
$i:=i+1$.
\end{itemize}

As a criterion of accuracy we can choose one of the following conditions:

1) the absolute error is less than $\Delta$;

2) the relative error is less than $\delta$;

3) both conditions 1) and 2) are satisfied.

Note that inequality (4) makes possible to check the absolute error without
a direct comparison with the initial number. This algorithm of rounding can
be applied to a result of any arithmetic operation if lengths of its
numerator and denominator exceed a given threshold  $M$. Of course, the
values of parameters $\Delta$, $\delta$, $M$ are given by the user.

{\bf 3. Estimations of round-off errors.} As a consequence of the
recurrence formula for the denominator the values $q_k$  are minimal for
each fixed $k$, if $a_i=1$ for $i=0,1,...,k$, i.e. $q_i=q_{i-1}+q_{i-2}$,
$q_{-2}=1$, $q_{-1}=1$, $q_0=q_{-1}+q_{-2}=1$. Therefore for any convergent
$\frac{p_k}{q_k}$ we obtain the inequality
$$
q_k \ge F_{k+1}, \eqno(5)
$$
where $F_k$ are Fibonacci numbers defined by the recurrence formulas
$F_k=F_{k-1}+F_{k-2}$ ($k\ge 2$), where $F_0=0$, $F_1=1$. It is well-known
that Fibonacci numbers are expressed by the formula:
$$
F_k=\frac{1}{\sqrt5}(\Phi^k-\hat\Phi^k), \eqno(6)
$$
where $\Phi=\frac{1}{2}(1+\sqrt5)\approx 1.618$ (``golden section''),
$\hat\Phi=\frac{1}{2}(1-\sqrt5)\approx -0.618$. It is also well-known that
the convergence of the continued fraction expansion of $\Phi$ is the
slowest among other numbers. From (5) and (6) it follows that
$$
q_kq_{k+1}\ge F_{k+1}F_{k+2}=\frac{1}{5}(\Phi^{k+1}-\hat\Phi^{k+1})
 (\Phi^{k+2}-\hat\Phi^{k+2})=
$$
$$
=\frac{1}{5}(\Phi^{2k+1}-(\Phi\hat\Phi)^{k+1}(\Phi+\hat\Phi)+\hat\Phi^{2k+1}) 
=\frac{1}{5}(\Phi^{2k+3}+(-1)^k+\hat\Phi^{2k+3}),
$$  
since $\Phi\hat\Phi=-1$ and $\Phi+\hat\Phi=1$. Using the values of $\Phi$ 
and $\hat\Phi$, we obtain for arbitrary $k$ the estimations 
$q_kq_{k+1}>\frac{1}{5}\Phi^{2k+2}$, while for even $k$: 
$q_kq_{k+1}>\frac{1}{5}\Phi^{2k+3}$. From these estimates and 
inequalities (4) it follows that the absolute round-off error is less 
than $\Delta$ if 
$$
k\ge \frac{1}{2}\log_\Phi\frac {5}{\Delta}-1 
$$
for even $k$ and 
$$
k\ge\frac{1}{2}\log_\Phi\frac{5}{\Delta}-\frac{3}{2}
$$
for odd $k$.

{\bf 4. Estimations of the number of iterations.} These relations lead to 
upper estimations for the number of iterations required for the 
approximation of a rational fraction with an absolute error smaller 
than $\Delta$. Let $\frac{p_k}{q_k}$ be the convergent of a continued 
fraction within the required error and the convergent 
$\frac{p_{k-1}}{q_{k-1}}$ does not give the required accuracy yet. 
Thus $k$ is the number of iterations necessary to obtain the required accuracy.

For any real number $r$ denote by $\lfloor r\rfloor$ the floor of $r$ and 
by $\lceil r\rceil$ the ceiling of $r$. So $\lceil r\rceil$ is the least
integer that is greater or equal to $r$; similarly, $\lfloor r\rfloor$ is
the integer part of $r$, i.e. the largest integer that is less or equal
to $r$. Hence if $r$ is an  integer number, then $\lfloor r \rfloor 
= \lceil r\rceil$ and otherwise $\lceil r\rceil= \lfloor r\rfloor+1
=\lfloor r+1\rfloor$. Then the required upper estimate has the following form:
$$
k\leq \lceil\frac{1}{2}\log_\Phi\frac
{5}{\Delta}-1\rceil\leq\lfloor\frac{1}{2}\log_\Phi\frac {5}{\Delta}
\rfloor. \eqno(7)
$$
But if the number $\lceil\frac{1}{2}\log_\Phi\frac {5}{\Delta}
-\frac{3}{2}\rceil$ is even, the estimation (7) can be strengthened:
$$
k\leq \lceil\frac{1}{2}\log_\Phi\frac
{5}{\Delta}-\frac{3}{2}\rceil\leq\lfloor\frac{1}{2}\log_\Phi
\frac {5}{\Delta}-\frac{1}{2}\rfloor \eqno(8)
$$
In the case of $\Delta=10^{-N}$ the estimations (7) and (8) give the 
following result:

{\bf Theorem.}
{\it If an absolute error specified in the criterion of accuracy of 
rounding has  the form } $\Delta=10^{-N}$, {\it then  }
$$
k\leq \lfloor a+bN\rfloor,\eqno(9)
$$
{\it where } $a=\frac{1}{2}\log_{\Phi}5\approx 1.672$ { \it and } 
$b=\frac{1}{2}\log_{\Phi}10\approx 2.392$. {\it If the number } 
$\lceil a+bN-\frac{3}{2}\rceil$ {\it is even, the estimation  } (9) 
{\it can be strengthened:}
$$
k\leq \lfloor a-\frac12 + bN\rfloor.\eqno(10)
$$

For example, if $N=8$, then the estimation (9) shows that $k\leq20$; 
for $N=9$ this estimation gives $k\leq 23$, but in this case the 
estimation (10) is applicable, so $k\leq 22$.

Note that these estimations depend only on the absolute error $\Delta$ 
and do not depend on the initial (i.e. rounded) numbers. In fact (see 
below), the number of iterations is usually much less than right-hand 
sides of these inequalities. Since the number of iterations is estimated 
by a linear function of logarithm of absolute error, an increase in the 
accuracy of calculations does not lead to an explosive growth of 
calculation time.

Heuristically, it is easy to estimate the mean {\it value} $\overline k$ 
of parameter $k$ for a fixed absolute error $\Delta$. Consider a convergent 
of a continued fraction $\frac{p_k}{q_k}$ as an approximation to a real 
number $x$.  A.~Ya.~Khinchin investigated the convergents of continued 
fraction $\frac{p_k}{q_k}$ for real numbers and proved that for  almost 
all $x$ the estimation $\lim_{k\to\infty}\sqrt[k]{q_k}=\gamma$ is valid, 
where $\gamma$ is a constant (see \cite{10}). P. Levy (\cite{9}, p. 320), 
showed that $\ln \gamma=\frac{\pi^2}{12\ln 2}\approx 1.18657\dots$, i.e. 
$\gamma\approx 3.27582\dots$ Roughly speaking, this result means that if 
values of $k$ are sufficiently large, then the denominator $q_k$ of a 
continued fraction  is ``close'' to $\gamma^k$. 

Leaving the mathematical rigor aside for a moment, substitute the 
quantities $\gamma^k$ and $\gamma^{k+1}$ into (4) for of  $q_k$  and 
$q_{k+1}$  in order to estimate a mean order of the convergent with a 
given approximation error $\Delta$. As an upper bound we obtain the 
number $\frac{\ln(1/\Delta)}{2\ln\gamma}-\frac{1}{2}$, and the lower 
bound differs from the upper bound by the value 
$\frac{\ln(1+1/\gamma)}{2 \ln\gamma}\approx 0.11$. Thus, the mean 
value of $k$ (not necessarily integer) is close to 
$\frac{\ln(1/\Delta)}{2\ln\gamma}$. If $\Delta=10^{-N}$, then
$$
\overline k\sim\frac{\ln(1/\Delta)}{2\ln\gamma}=\frac{N\ln 10}{2\ln\gamma} 
\approx 0,97\cdot N\sim N.\eqno(11)
$$
This estimation becomes realistic only for large values of $N$, otherwise 
$\overline k$ is much less than $N$.

{\bf 5. Examples of applications of different variants of approximate 
rational arithmetic.}  To compare different variants of rational arithmetics
consider a classical example of a numerical calculation of the function 
$\sin x$ at points $x_m=\frac{\pi}{6}+2\pi m$ by summation of its Taylor 
series. The sum is calculated until the absolute value of a summand 
becomes less than $10^{-7}$. The number $\pi$ is replaced by its rational 
approximation $\frac{355}{113}$ with an absolute error $2.7 \cdot 10^{-7}$.

A 500 MHz Intel Pentium III processor  was used for calculations.
Different variants of rational arithmetics were created using the 
arbitrary precision arithmetic, implemented by means of the $C^{++}$ 
programming language  (implementation by means of the REDUCE  
system gives about the similar results).

We consider the following variants of approximate rational arithmetic:

I) The arbitrary precision arithmetic (without rounding).

II) Approximate rational arithmetic described in the section 2 with 
$M=9$, $\Delta=10^{-8}$, $\delta=\infty$ (so only absolute round-off 
error $\Delta=10^{-8}$ is fixed).

III) The same arithmetic with $M=9$, $\Delta=\delta=10^{-8}$.

IV) The same arithmetic with $M=9$, $\Delta=\infty$, $\delta=10^{-8}$ 
(so only relative round-off error is fixed).

V) Fixed slash arithmetic \cite{2, 3}, where the maximum length $L$ of 
numerator and denominator is fixed by $L=6$.

\noindent
\begin{tabular}{|c|l|c|c|c|c|c|c|c|}
\hline
& m & 0 & 1 & 2 & 3 & 4 & 5 & 6  \\
\hline
I &$\varepsilon$ &$4\cdot 10^{-8}$&$4\cdot 10^{-7}$&$8\cdot10^{-7}$&$10^{-6}$&$2\cdot 10^{-6}$&$3\cdot 10^{-6}$&$5\cdot 10^{-6}$\\
     & s                 &$62$&$214$&$372$&$504$&$650$&$810$&$980$\\
     & t & $0.007$&$0.09$&$0.24$&$0.95$&$3.3$&5.7&$17$\\
\hline

II &$\varepsilon$&$2\cdot 10^{-8}$&$5\cdot 10^{-7}$&$10^{-6}$&$10^{-6}$&$2\cdot 10^{-6}$&$2\cdot 10^{-
6}$&$3\cdot 10^{-6}$ \\
$\Delta=10^{-8}$     & s                 & $16$&$13$&$12$&$12$&$12$&$12$&$11$\\
     & t & $0.007$&$0.08$&$0.15$&$0.21$&$0.28$&$0.34$&$0.42$\\
\hline
III &$\varepsilon$&$4\cdot 10^{-8}$&$5\cdot 10^{-7}$&$10^{-6}$&$10^{-6}$&$2\cdot 10^{-6}$&$2\cdot 
10^{-6}$&$3\cdot 10^{-6}$ \\
$\Delta=10^{-8}$      & s                 &15&13&12&$12$&$12$&$12$&$11$\\
$\delta=10^{-8}$      & t &$0.012$&$0.09$&$0.16$&$0.23$&$0.32$&$0.37$&$0.46$\\
\hline
IV &$\varepsilon$&$4\cdot 10^{-8}$&$3\cdot 10^{-7}$&$3\cdot 10^{-4}$&$0.21$&$0.6$&$0.8$&$1.17$\\
$\delta=10^{-8}$     & s                  &15&13&9&$9$&$9$&$8$&$8$\\
    & t & $0.014$&$0.06$&$0.14$&$0.18$&$0.25$&$0.29$&$0.34$\\
\hline
V &$\varepsilon$&$0$&$0$&$10^{-3}$&$0.7$&$1.0$&$1.4$&$3.4$ \\
L=6     & s                 & $2$&$2$&$12$&$12$&$12$&$11$&$12$\\
     & t & $0.017$&$0.12$&$0.18$&$0.19$&$0.21$&$0.23$&$0.28$\\
\hline
VI &$\varepsilon$&$4\cdot 10^{-8}$&$5\cdot 10^{-7}$&$10^{-6}$&$0.008$&$0.08$&$0.3$&$0.6$\\
L=9     & s                 & $17$&$18$&$17$&$18$&$18$&$16$&$18$\\
     & t & $0.025$&$0.21$&$0.29$&$0.31$&$0.33$&$0.36$&$0.39$\\
\hline
VII &$\varepsilon$&$4\cdot 10^{-8}$&$5\cdot 10^{-7}$&$10^{-6}$&$10^{-6}$&$10^{-6}$&$2\cdot 10^{-
4}$&$0.007$\\
L=12     & s                 & $24$&$24$&$23$&$23$&$24$&$23$&$24$\\
     & t & $0.05$&$0.29$&$0.49$&$0.56$&$0.64$&$0.65$&$0.67$\\
\hline
VIII &$\varepsilon$&$4\cdot 10^{-8}$&$5\cdot 10^{-7}$&$2\cdot10^{-6}$&$10^{-4}$&$0.04$
&$0.06$&$0.8$\\
S=12     & s                 & $11$&$11$&$10$&$11$&$11$&$10$&$11$\\
     & t & $0.013$&$0.18$&$0.34$&$0.48$&$0.51$&$0.52$&$0.54$\\
\hline
IX &$\varepsilon$&$4\cdot 10^{-8}$&$5\cdot 10^{-7}$&$10^{-6}$&$6\cdot 10^{-6}$&$2\cdot 10^{-
3}$&$0.01$&$0.4$\\
S=15     & s                 & $14$&$14$&$13$&$13$&$14$&$13$&$13$\\
     & t & $0.05$&$0.21$&$0.37$&$0.41$&$0.56$&$0.61$&$0.64$\\
\hline
X &$\varepsilon$&$4\cdot 10^{-8}$&$5\cdot 10^{-7}$&$10^{-6}$&$2\cdot 10^{-6}$&$3\cdot 10^{-6}$&$3\cdot 10^{-5}$&$0.01$\\
S=18     & s                 & $17$&$17$&$15$&$17$&$16$&$17$&$17$\\
     & t &$0.023$&$0.31$&$0.51$&$0.68$&$0.75$&$0.85$&$0.87$\\
\hline
XI &$\varepsilon$&$0$&$4\cdot 10^{-5}$&$0.04$&$0.1$&$0.2$&$0.4$&$0.9$\\
D=9     & s                  &$18$&$18$&$18$&$18$&$18$&$18$&$18$\\
     & t &$0.008$&$0.04$&$0.12$&$0.17$&$0.26$&$0.33$&$0.42$\\
\hline
\end{tabular}

\bigskip

\centerline
{\bf Table 1.}

\centerline
{\bf Comparison of rational arithmetics}

\bigskip

VI) The same arithmetic with $L=9$.

VII) The same arithmetic with $L=12$.

VIII) Floating-slash arithmetic \cite{2, 3}, where the maximum sum $S$ of  
lengths of numerator and denominator is fixed by $S=12$.

IX) The same arithmetic with $S=15$.

X) The same arithmetic with $S=18$.

XI) A ``reductive'' arithmetic, where the numbers of true top digits $D$ 
of numerator and denominator are fixed by $D=9$ (the other digits are 
replaced by zeros).

The results of numerical computations are presented in Table 1.
For calculation of the function $\sin(\frac{\pi}{6}+2\pi m)$ at points 
$m=0,1,2,3,4,5,6$ an absolute errors of the result $\varepsilon$ (the 
relative error equals $2\varepsilon$), time of calculations $t$ in seconds 
and the sum of the lengths of numerator and denominator $s$ for the result 
are specified in Table 1. All values of errors are rounded to the first 
digits to fit them in the format of the table.

It follows from Table 1 that in the case of infinite precision arithmetic 
the quick increase of the parameter $s$ leads to the explosive growth of 
calculations time. Curiously enough, the arbitrary precision arithmetic 
does not always give the most accurate result. This phenomenon can be 
partly explained by the inaccuracy of presentation of the number $\pi$, 
but it is also a consequence of the unusually simple form of the rational 
fraction $\sin(\frac{\pi}{6}+2\pi n)=\frac{1}{2}$. The nature of this 
effect is discussed in \cite{2}. If only the relative error of rounding is 
fixed (variant IV), the situation characteristic for the floating point 
arithmetic (variant XI) is repeated completely: starting at $m=4$ the 
errors are larger then half of computed values (see. \cite{10}, part 3).

For other types of rational arithmetics the increase of the difficulty
of calculations damages the accuracy of the result essentially,
in contrast to variants II and III of approximate rational arithmetic. 
On the other hand, the time of calculations is comparable in all cases 
except of the arbitrary precision arithmetic.

\bigskip
\noindent
\begin{tabular}{|c|c|c|c|c|c|c|c|c|c|c|c|}
\hline
 N &  16 & 18 & 20 & 22 & 24 & 26 & 28 & 30 & 32 & 34 & 36 \\
\hline
 &&&&&&&&&&& \\
$\bar k$ &13.9&16.4&18.9&20.9&22.7&24.7&26.5&28.4&30.9&33.0&34.9\\
\hline
\end{tabular}

\bigskip

\centerline
{\bf Table 2. Dependence of the mean number }
\centerline
{\bf of iterations on the accuracy of rounding}

\bigskip

Table 2 presents the dependence of the mean number of iterations $\bar k$ 
(see above, Section 4) on the absolute rounding error $\Delta=10^{-N}$ for
$N=16,18,...,36$. The calculations were implemented for the model example 
described above with $M=9$. It follows from Table 2 that an estimation (11) 
is realistic for $\bar k$ if $N$ is sufficiently large.

{\bf 6. Conclusion.} The approximate rational arithmetic described in the 
section 2 provides a sufficiently higher degree of accuracy of calculations 
in a time comparable to other types of rational arithmetics. This result is 
illustrated by the above calculations quite clearly. It is particularly 
important that the round-off error can be controlled by the user on each 
step of the calculation procedure. This allows us to control the inaccuracy 
of rounding, estimate the maximum computing error beforehand, and guarantee 
(in particular, in terms of Interval Analysis) the required accuracy of 
calculations.

\end{document}